\documentclass[12pt,reqno]{amsart}
\usepackage{amsthm,amsfonts,amssymb,amscd}
\newtheorem{theorem}[subsection]{Theorem}
\newtheorem*{theorema}{Theorem A}
\newtheorem{proposition}[subsection]{Proposition}
\newtheorem*{proposition*}{Proposition}
\newtheorem{lemma}[subsection]{Lemma}
\newtheorem{corollary}[subsection]{Corollary}
\theoremstyle{definition}
\newtheorem{definition}[subsection]{Definition}
\newtheorem{example}[subsection]{Example}
\theoremstyle{remark}
\newtheorem{remark}[subsection]{Remark}
\newtheorem*{remark*}{Remark}
\newcommand{\mt}[1]{\operatorname{#1}}

\newcommand{\NE}{{NE}}
\newcommand{\Diff}{\operatorname{Diff}}
\newcommand{\Supp}{\operatorname{Supp}}

\newcommand{\Const}{\operatorname{Const}}
\newcommand{\qq}{\mathbin{\sim_{\scriptscriptstyle{\QQ}}}}
\newcommand{\pal}{\text{---}}

\newcommand{\bir}{\dasharrow}
\newcommand{\ov}[1]{\overline{#1}}
\newcommand{\down}[1]{\lfloor #1\rfloor}
\newcommand{\up}[1]{\lceil #1\rceil}
\newcommand{\fr}[1]{\{ #1\}}
\renewcommand{\AA}{{\mathbb A}}
\newcommand{\CC}{{\mathbb C}}
\newcommand{\RR}{{\mathbb R}}

\newcommand{\QQ}{{\mathbb Q}}
\newcommand{\PP}{{\mathbb P}}
\newcommand{\NN}{{\mathbb N}}

\newcommand{\DD}{{\mathbb D}}
\newcommand{\EE}{{\mathbb E}}
\newcommand{\OOO}{{\mathcal O}}
\newcommand{\AAA}{{\mathcal A}}
\newcommand{\MMM}{{\mathcal M}}

\newcommand{\CCC}{{\mathcal C}}

\newcommand{\NNN}{{\mathcal N}}
\newcommand{\ep}{\varepsilon}
\newcommand{\var}{\varphi}

\date{}
\title{Boundedness of non-birational extremal contractions}
\author{Yuri~G.~Prokhorov\thanks{Partially supported
by the Russian Foundation of Fundamental Research}}

\email{prokhoro@mech.math.msu.su} 
\address{Department of Mathematics
(Algebra Section), Moscow State University, 117234 Moscow, Russia}

\address{Current address:\quad Department of Mathematics,
Tokyo Institute of Technology, Oh-okayama, Meguro, 152 Tokyo,
Japan}

\begin{document}
\begin{abstract}
We consider $K_X$-negative extremal contractions \mbox{$f\colon
X\to (Z,o)$}, where $X$ is an algebraic threefold with only
$\epsilon$-log terminal $\QQ$-factorial singularities and $(Z,o)$
is a two (resp., one)-dimensional germ. The main result is that
$K_X$ is $1$, $2$, $3$, $4$ or $6$-complementary or we have, so
called, exceptional case and then the singularity $(Z\in o)$ is
bounded (resp., the multiplicity of the central fiber $f^{-1}(o)$
is bounded).
\end{abstract}
\maketitle

\section*{Introduction}
The aim of this paper is to generalize in some sense the following
two statements connected with the \textit{exceptionality}
phenomenon which was discovered by Shokurov~\cite{Sh}. For the
definition of complements we refer to~\ref{def-compl}.

\begin{enumerate}
\item[1)]
The class of two-dimensional $\ep$-log terminal singularities
which have no $1$ or $2$-complements is bounded~\cite{Sh}, see
also~\cite{MP}.
\item[2)]
Numerical invariants of any three-dimensional $\ep$-log terminal
singularity which have no $1$, $2$, $3$, $4$ or $6$-complements
are bounded~\cite{Sh1}, see also~\cite{MP} for the case of
quotient singularities.
\end{enumerate}
There are generalizations of these facts for cases of
two-dimensional elliptic fibrations (see~\cite{Sh1}) and
three-dimensional birational contractions~\cite{Sh1}. In this
paper we apply Shokurov's inductive approach to the study of
contractions of relative dimension one or two. Such contractions
naturally appear at the end of log MMP for varieties of Kodaira
dimension $\kappa=-\infty$ (see~\cite{KMM}). Our main result is
the following

\begin{theorema}
Let $X$ be an algebraic threefold with only $\ep$-log terminal
${\QQ}$-factorial singularities and let $f\colon X\to Z$ be a
$K_X$-negative extremal contraction (that is a projective morphism
such that $f_*\OOO_X=\OOO_Z$, $\rho (X/Z)=1$ and $-K_X$ is
$f$-ample). Assume that $\dim Z=2$ or $1$. Fix a point $o\in Z$.
Then one of the following holds.
\begin{enumerate}
\item
There exists an $1$, $2$, $3$, $4$ or $6$-complement of $K_X$ near
$f^{-1}(o)$. Moreover, there exists such a complement which is not
plt and even non-exceptional in the sense of \ref{def-sing}.
\item
$K_X$ is exceptionally $n$-complementary near $f^{-1}(o)$ for
$n<\NNN_2$, where $\NNN_2$ is a constant which does not depend on
$f$. Then we have
\begin{enumerate}
\item
in the case $\dim Z=2$ the singularity $(Z\in o)$ is bounded,
i.~e. up to analytic isomorphisms it belongs to a finite set
$\MMM(\ep)$ depending only on $\ep$;
\item
in the case $\dim Z=1$ the multiplicity of the central fiber
$f^{-1}(o)$ is bounded by a constant $\AAA(\ep)$ depending only on
$\ep$.
\end{enumerate}
\end{enumerate}
\end{theorema}
For the first time the notion of complements was introduced by
Shokurov in \cite{Sh}. Roughly speaking an $n$-complement of $K_X$
is a ``good'' divisor in the multiple anticanonical linear system
$|-nK_X|$, see \ref{def-compl}. The set of numbers $\{1, 2, 3, 4,
6\}$ from (i) of Theorem~A is well known in the theory of
algebraic surfaces, these numbers are called \textit{regular}.

Following Shokurov we call case (ii) of Theorem~A
\textit{exceptional}. More precisely, a contraction $f\colon X\to
Z$ is said to be \textit{exceptional} near $o\in Z$ if for any
complement $K_X+D$ on $X$ which is not Kawamata log terminal,
there exists exactly one divisor $S$ of $K(X)$ with discrepancy
$a(S,D)=-1$. In this case we also show in Proposition~\ref{unique}
that the divisor $S$ does not depend on the choice of $D$, so it
is, in some sense, a distinguished divisor of $K(X)$. Thus
Theorem~A states that either $K_X$ has a regular non-exceptional
complement (see \ref{def-sing}) near the fiber over $o\in Z$ or
$f\colon X\to Z$ is exceptional and then we have case (ii).
Shokurov in~\cite[\S 7]{Sh1} investigated from this point of view
the case when $f$ is a birational contraction.

There is a two-dimensional analog of Theorem~A:

\begin{proposition*}[\cite{Pr}]
Let $X$ be an algebraic surface having only quotient singularities
and let $f\colon X\to Z$ be a $K_X$-negative extremal contraction
onto a curve. Fix a point $o\in Z$. Then one of the following
holds.
\begin{enumerate}
\item
There exists a non-exceptional $1$ or $2$-complement of $K_X$ near
$f^{-1}(o)$.
\item
$K_X$ is exceptionally $3$-complementary near $f^{-1}(o)$. In this
case $X$ has exactly two singular points on $f^{-1}(o)$ which are
of types $\frac{1}{2b-3}(1,b-2)$, $b\ge 2$ and
\begin{equation*}
\begin{array}{ccccccc}
&&&&\stackrel{-3}{\circ}&&\\ &&&&|&&\\
\stackrel{-2}{\circ}&\pal&\stackrel{-2}{\circ}&\pal&\stackrel{-b}{\circ}
&\pal&\stackrel{-2}{\circ}\\
\end{array}
\end{equation*}
\end{enumerate}
\end{proposition*}
We hope that similar to the proposition above the exceptional
cases in (ii) of Theorem~A can be classified.

The most important case for applications is when $\ep=1$, i.~e.
when $X$ has only terminal singularities. In this case we expect
that $K_X$ is $1$-complementary (this is a weak form of Reid's
``general elephant'' conjecture). Assuming that this conjecture is
true in the case $\dim Z=2$ we have a rough classification of
contractions $f\colon X\to Z$~\cite{Pr0}. In particular, the base
surface $S$ has only Du Val singularities of type $A_n$
(Iskovskikh's conjecture,~\cite{I}, \cite{Pr0}). These conjectures
have applications to the rationality problem of conic bundles
(see~\cite{I}).

The idea of the proof of Theorem~A is an application of Shokurov's
theorem on boundedness of two-dimensional complements. We
construct some special, so called plt blow-up of $X$
(see~\cite{Pr1}) and show that, in some sense, the geometry of
exceptional divisor reflects the geometry of $X$ itself.

The paper is organized as follows: Sect.~1 is auxiliary. In
Sect.~2 we construct a new birational model of $X$ so that we
reduce the problem to dimension $2$. In Sect.~3 we prove
Theorem~A. Examples in Sect.~\ref{last} show that boundedness as
in Theorem~A under weaker restrictions cannot be expected. Some of
the results were announced in~\cite{Pr2}.

All varieties are assumed to be algebraic and defined over $\CC$,
the field of complex numbers. A \textit{contraction} (or
\textit{extraction}, if we start with $X$ instead of $Y$) is a
projective morphism of normal varieties $f\colon Y\to X$ such that
$f_*\OOO_Y=\OOO_X$. A \textit{blow-up} is a birational extraction.

\section{Preliminary results}

All the necessary facts and definitions from log MMP can be found
in~\cite{KMM}, \cite{Sh}, \cite{Ut} and \cite{Sh2}. We follow
essentially the terminology and notation of \cite{Ut}, \cite{Sh}
and \cite{Sh2} (see also~\cite{Ko} for a nice introduction to
singularities of pairs).

\begin{definition}\label{def-sing}
Let $X$ be a normal algebraic variety and let $D=\sum d_iD_i$ be a
$\QQ$-divisor on $X$. $D$ is called a \textit{boundary} (resp.
\textit{subboundary}) if $0\leq d_i\leq 1$ (resp. $d_i\leq 1$) for
all $i$. Let $f\colon Y\to X$ be a projective birational morphism.
Assume that $K_X+D$ is $\QQ$-Cartier and write
\begin{equation*}
\label{form-def}
K_Y\equiv f^*(K_X+D)+\sum_E a(E,D)E,
\end{equation*}
where $E$ runs over prime divisors on $Y$, $a(E,D)\in\QQ$, and
$a(D_i,D)=-d_i$ for each component $D_i$ of $D$. The coefficients
$a(E,D)$ are called \textit{discrepancies} of $(X,D)$. Define
\begin{equation*}
\mt{discr}(X,D):=\inf_{E,f}\{a(E,D): E\ \text{is $f$-exceptional
divisor} \}.
\end{equation*}
Then the pair $(X,D)$ or, by abuse of language, the divisor
$K_X+D$ is said to be\par
\begin{tabular}{ll}
\textit{Kawamata log terminal (klt)} &iff $\mt{discr}(X,D)>-1$ and
$\down{D}\le 0$,\\ \textit{purely log terminal (plt)} &iff
$\mt{discr}(X,D)>-1$,\\ \textit{$\ep$-log terminal ($\ep$-lt)}
&iff $\mt{discr}(X,D)>-1+\ep$,\\ \textit{log canonical (lc)} &iff
$\mt{discr}(X,D)\ge -1$.\\
\end{tabular}
\par
A pair $(X,D)$ is said to be \textit{divisorial log terminal
(dlt)} if there exists a good resolution $f\colon Y\to X$ such
that the exceptional locus consists of divisors with $a(E,D)>-1$.

A log canonical pair $(X,D)$ is said to be \textit{exceptional} if
there exists \textit{exactly one} divisor $E$ with discrepancy
$a(E,D)=-1$. So if $(X,D)$ is exceptional, then either $(X,D)$ is
plt and $\down{D}=E$ or $\down{D}=\emptyset$ and $E$ is an unique
exceptional divisor with $a(E,D)=-1$.
\end{definition}

Everywhere below if we do not specify the opposite we consider log
pairs $(X,D)$ consisting of a normal variety $X$ and a boundary
$D$ on it.

\begin{definition}[{\cite[\S 3]{Sh}}, {\cite[Ch. 16]{Ut}}]
Let $X$ be a normal variety, $S\ne\emptyset$ be an effective
reduced divisor on $X$, and let $B$ be a $\QQ$-divisor on $X$,
such that $S$, $B$ have no common components. Assume that $K_X+S$
is lc in codimension two. Then the \textit{different} of $B$ on
$S$ is defined by
\begin{equation*}
K_S+\Diff_S(B)\equiv (K_X+S+B)|_S.
\end{equation*}
Usually we will write simply $\Diff_S$ instead of $\Diff_S(0)$
\end{definition}

\begin{theorem}[Inversion of Adjunction
{\cite[3.3]{Sh}}, {\cite[17.6]{Ut}}]
\label{Inv-Adj}
Let $X$ be a normal variety and let $D$ be a boundary on it. Write
$D=S+B$, where $S=\down{D}$ and $B=\fr{D}$. Assume that $K_X+S+B$
is $\QQ$-Cartier. Then $(X,S+B)$ is plt near $S$ iff $S$ is normal
and $(S,\Diff_S(B))$ is klt.
\end{theorem}

\begin{definition}[\cite{Pr1}]
Let $X$ be a normal variety and let $g\colon Y\to X$ be a blow-up
such that the exceptional locus of $g$ contains only one
irreducible divisor, say $S$. Assume that $K_Y+S$ is plt and
$-(K_Y+S)$ is $f$-ample. Then $g\colon (Y\supset S)\to X$ is
called a \textit{purely log terminal (plt) blow-up} of $X$.
\end{definition}

\begin{definition}[{\cite[3.14]{Sh}}]
Let $X$ be a normal variety and let $D=\sum d_iD_i$ be a
$\QQ$-divisor on $X$ such that $K_X+D$ is $\QQ$-Cartier. A
subvariety $W\subset X$ is said to be a \textit{center of log
canonical singularities} for $(X,D)$ if there exists a prime (not
necessary exceptional) divisor $E$ over $X$ with center at $W$ and
discrepancy $a(E,D)\le-1$. The union of all centers of log
canonical singularities is called the \textit{locus of log
canonical singularities} and denoted by $LCS(X,D)$.
\end{definition}

For the following statement we refer to~\cite[17.4]{Ut} (in
dimension $2$ it was proved earlier by Shokurov~\cite{Sh}).

\begin{theorem}[Connectedness Lemma]
\label{connect}
Let $X$ be a normal projective variety, let $f\colon X\to Z$ be a
contraction and let $D=\sum d_iD_i$ be a boundary on $X$ such that
$K_X+D$ is $\QQ$-Cartier. If $-(K_X+D)$ is $f$-nef and $f$-big,
then $LCS(X,D)$ is connected near each fiber of $f$.
\end{theorem}

V.~V.~Shokurov informed me that the result above has a
generalization modulo log MMP to the case when $-(K_X+D)$ is only
nef (preprint, in preparation,
cf.~\cite[6.9]{Sh}\footnote{Recently a generalization of this fact
was obtained by O.~Fujino \cite[Proposition~2.1]{F}}). We need
only a particular case of this fact (see \ref{connect-(2)3}).

\begin{definition}[{\cite[5.1]{Sh}}]
\label{def-compl}
Let $X$ be a normal variety and let $D=S+B$ be a subboundary, such
that $B$, $S$ have no common components, $S$ is a reduced divisor,
and $\down{B}\le 0$. Then one says that $K_X+D$ is
\textit{$n$-complementary}, if there exists a $\QQ$-divisor $D^+$
such that
\begin{enumerate}
\item
$nD^+$ is integer and $n(K_X+D^+)\sim 0$;
\item
$K_X+D^+$ is lc;
\item
$nD^+\geq nS+\down{(n+1)B}$.
\end{enumerate}
The divisor $K_X+D^+$ is called an \textit{$n$-complement} of
$K_X+D$.
\end{definition}

\begin{remark}
Note that in general it is not true that $D^+\ge D$. This however
is true if all the coefficients of $D$ are \textit{standard},
i.~e. they have the form $d_i=1-1/m_i$, where
$m_i\in\NN\cup\{\infty\}$~\cite[2.7]{Sh1}.
\end{remark}

For convenience, we recall several facts about complements. For
proofs we refer to \cite{Sh}, \cite{Ut}, \cite{Sh1} and author's
lecture notes \cite{Pr-l}.

\begin{lemma}[{\cite[5.4]{Sh}}]\label{down}
Let $f\colon Y\to X$ be a blow-up and let $D$ be a subboundary on
$Y$. Assume that $K_Y+D$ is $n$-complementary. Then $K_X+f_*D$ is
$n$-complementary.
\end{lemma}

\begin{lemma}[{\cite[4.4]{Sh1}}]
\label{podnyal}
Let $f\colon Y\to X$ be a blow-up and let $D$ be a subboundary on
$Y$ such that
\begin{enumerate}
\item
$K_Y+D$ is $f$-nef;
\item
$f_*D=\sum d_if_*D_i$ is a boundary with standard coefficients.
\end{enumerate}
Assume that $K_X+f_*D$ is $n$-complementary and let $K_X+(f_*D)^+$
be any $n$-complement. Then its crepant pull-back $K_X+(f_*D)^+$
is an $n$-complement of $K_Y+D$.
\end{lemma}

\begin{lemma}[{\cite[2.1]{Pr}, cf.~\cite[19.6]{Ut}}]
\label{prodolj}
Let $(X,S)$ be a plt pair with reduced $S\ne 0$. Let $f\colon X\to
Z$ be a contraction such that $-(K_X+S)$ is $f$-nef and $f$-big.
Fix a fiber $f^{-1}(o)$, $o\in Z$ meeting $S$. Assume that
$K_S+\Diff_S$ is $n$-complementary. Then any $n$-complement of
$K_S+\Diff_S$ in a neighborhood of $f^{-1}(o)$ can be extended to
an $n$-complement of $K_X+S$. This means that for any
$n$-complement $K_S+\Diff_S^+$ there exists an $n$-complement
$K_X+S+B^+$ such that $\Diff_S(B^+)=\Diff_S^+$.
\end{lemma}

\begin{proof}
Let $K_S+\Theta$ be an $n$-complement of $K_S+\Diff_S$. Consider a
good resolution $g\colon Y\to X$ and put $h=f\circ g$. Write
\begin{equation*}
g^*(K_X+S)=K_Y+S_Y+B,
\end{equation*}
where $S_Y$ is the proper transform of $S$ and $B$ is a
subboundary on $Y$ such that $\down{B}\le 0$. Then $K_Y+S_Y+B$ is
plt \cite[3.10]{Ko}. By Inversion of Adjunction
$g^*(K_S+\Diff_S)=K_{S_Y}+\Diff_{S_Y}(B)$ is klt and since $Y$ is
non-singular, $\Diff_{S_Y}(B)=B|_{S_Y}$. By Lemma~\ref{podnyal}
$K_{S_Y}+\Theta_Y:=g^*(K_S+\Theta)$ is an $n$-complement of
$K_{S_Y}+\Diff_{S_Y}(B)$. In particular,
$\ov{\Theta_Y}:=n\Theta_Y-\down{(n+1)\Diff_{S_Y}(B)}$ is an
effective divisor from the linear system
$|-nK_{S_Y}-\down{(n+1)B}|_{S_Y}|$ (see (iii) of \ref{def-compl}).

From the exact sequence
\begin{multline*}
0\longrightarrow\OOO_{Y}(-nK_Y-(n+1) S_Y-\down{(n+1)B})\\
\longrightarrow\OOO_{Y}(-nK_Y-nS_Y-\down{(n+1)B})\\
\longrightarrow\OOO_{S_Y}(-nK_{S_Y}-\down{(n+1)B}|_{S_Y})
\longrightarrow 0
\end{multline*}
and the vanishing
\begin{multline*}
R^1h_*(\OOO_{Y}(-nK_Y-(n+1)S_Y-\down{(n+1) B}))=\\
R^1h_*(\OOO_{Y}(K_Y+\up{-(n+1)(K_Y+S_Y+B)}))=0.
\end{multline*}
We get surjectivity of the restriction
\begin{multline*}
H^0(Y,\OOO_{Y}(-nK_Y-nS_Y-\down{(n+1)B})) \longrightarrow\\
H^0(S_Y,\OOO_{S_Y}(-nK_{S_Y}-\down{(n+1)B}|_{S_Y})).
\end{multline*}
Therefore there is an element $\ov{D}\in |-nK_Y-nS_Y-\down{(n+1)
B}|$ such that $\ov{D}|_{S_Y}=\ov{\Theta_Y}$. Put
$B_Y^+:=\frac{1}{n}(\down{(n+1)B}+\ov{D})$ and $B^+:=g_*B_Y^+$.
Then we have $n(K_Y+S_Y+B_Y^+)\sim 0$ and
$(K_Y+S_Y+B_Y^+)|_{S_Y}=K_{S_Y}+\Theta_Y$. This gives us
$n(K_X+S+B^+)\sim 0$ and $(K_X+S+B^+)|_S=K_S+\Theta$. By our
construction $B^+$ is effective, so it is a boundary. Now we have
to show only that $K_X+S+B^+$ is lc. Since $K_S+\Diff_S$ is klt
and $K_S+\Theta$ is lc, we have that
$K_S+\Diff_S+\alpha(\Theta-\Diff_S=(K_X+S+\alpha B^+)|_S$ is klt
for $\alpha<1$. By Inversion of Adjunction $K_X+S+\alpha B^+$ is
plt near $S$ for $\alpha<1$. Therefore $K_X+S+B^+$ is lc near $S$.
Moreover, by Connectedness Lemma $LCS(X,S+\alpha B^+)$ is
connected near each fiber of $f$ for $0<\alpha<1$ (because
$-(K_X+S+\alpha B^+)$ is $f$-nef and $f$-big). Therefore
$LCS(X,S+\alpha B^+)=S$ and $K_X+S+\alpha B^+$ is plt for
$\alpha<1$. This gives us that $K_X+S+B^+$ is lc and proves the
lemma.
\end{proof}

\begin{theorem}[{\cite{Sh1}}]
\label{Const}
Let $S$ be a projective surface and let $\Delta$ be a boundary on
it. Assume that
\begin{enumerate}
\item
coefficients of $\Delta$ are \textsl{standard}, i.~e. $\Delta
=\sum (1-1/m_i)\Delta_i$, $m_i\in\NN\cup\{\infty\}$;
\item
$K_S+\Delta$ is lc;
\item
$-(K_S+\Delta)$ is nef and big.
\end{enumerate}
Then there exists an $n$-complement of $K_S+\Delta$ such that
$n<\NNN_2$, where $\NNN_2$ is an absolute constant which does not
depend on $(S, \Delta)$.
\end{theorem}

\section{Construction of a good model}
In this section we modify the techniques developed in \cite[\S
6]{Sh}, \cite[\S 21]{Ut} to our situation. Roughly speaking we
have to construct a ``good'' model for a plt blow-up of our
threefold $X$ from Theorem~A in order to apply Lemma~\ref{prodolj}
and Theorem~\ref{Const}. A similar techniques was applied by
Shokurov in~\cite[\S 7]{Sh1} for the study of the birational case.
However we employ here a little different method which allow us to
work with irreducible exceptional divisor. This construction gives
us the existence of an $n$-complement on $X$ which is exceptional
and $n<\NNN_2$.

\subsection{Notation}\label{notations}
Let $X$ be an algebraic threefold with only klt $\QQ$-factorial
singularities, $f\colon X\to Z$ be an extremal contraction on a
(normal) variety $Z$ of positive dimension, i.~e. we assume that
$\rho(X/Z)=1 $ and $-K_X$ is $f$-ample. Fix a point $o\in Z$ and
assume that $Z$ is a germ near $o$. In this section we do not
assume that $\dim Z\ne 3$. Thus we consider the following cases:
1) del Pezzo fibrations, 2) generically conic bundle fibrations,
3) birational contractions.

\subsection{}\label{notation1}
Take a boundary $F$ on $X$ such that
\begin{enumerate}
\item
$-(K_X+F)$ is $f$-ample;
\item
$K_X+F$ is lc but not klt near $f^{-1}(o)$.
\end{enumerate}

By our assumptions, $\dim Z\geq 1$, so we can take simply
$F=cf^*L$, where $L$ is an effective divisor on $Z$ containing $o$
and $c$ is the \textit{log canonical threshold} of $(X,f^*L)$ that
is maximal $c\in\QQ$ such that $K_X+cf^*L$ is lc. Of course there
are another possibilities for the choice of $F$.
\par
First we consider the easy case, when $K_X+F$ is plt. It is clear
that $\down{F}\neq\emptyset$, because $K_X+F$ is not klt. If $\dim
Z=2$ and under the assumption that $X$ has only terminal
singularities this case was studied in~\cite{Pr}.

\begin{lemma}[{cf.~\cite[5.12]{Sh}}]\label{semistable}
Notation as above. Assume that $K_X+F$ is plt (near $f^{-1}(o)$).
Then
\begin{enumerate}
\item
$K_X+\down{F}$ is $n$-complementary for $n<\NNN_2$;
\item
if $\down{F}$ contains a non-compact component, then
$K_X+\down{F}$ has a regular complement which is not plt.
\end{enumerate}
\end{lemma}
Note that by Connectedness Lemma~\ref{connect} $\down{F}$ is
connected near $f^{-1}(o)$, so it is irreducible. Moreover, if
$\down{F}$ is compact, then $\down{F}=f^{-1}(o)$. This is
impossible if $\dim Z=2$, because $f$ is equidimensional.

\begin{proof}
First we assume that $-(K_X+\down{F})$ is $f$-ample and so is
$-(K_{\down{F}}+\Diff_{\down{F}})$. By Adjunction~\ref{Inv-Adj}
$K_{\down{F}}+\Diff_{\down{F}}$ is klt and $\Diff_{\down{F}}$ has
only standard coefficients~\cite[3.9]{Sh}, \cite[16.6]{Ut}. By
Theorem~\ref{Const} $K_{\down{F}}+\Diff_{\down{F}}$ is
$n$-complementary for some $n<\NNN_2$. Moreover, we can take
$n\in\{1,2,3,4,6\}$ if ${\down{F}}$ is non-compact,
see~\cite[5.6]{Sh} for $\dim f({\down{F}})=2$ and~\cite[\S 3]{Sh1}
for $\dim f({\down{F}})=1$ (see also \cite{Ut} and \cite{Pr-l}).
Since $-(K_X+{\down{F}})$ is $f$-ample, by Lemma~\ref{prodolj}
complements of $K_{\down{F}}+\Diff_{\down{F}}$ can be extended to
complements of $K_X+{\down{F}}$ near the fiber over $o\in Z$.
\par
Now we assume that $-(K_X+{\down{F}})$ is not $f$-ample. Then
$K_X+{\down{F}}$ is $f$-nef, because $\rho (X/Z)=1$. Therefore
$-\fr{F}$ is $f$-ample. This is possible only if $f$ is
birational. If $f$ is divisorial, then $Z$ is $\QQ$-factorial and
$K_Z+f_*F$ is plt, hence so is $K_Z+f_*{\down{F}}$. Moreover,
$\dim f({\down{F}})=2$, because $f$ contracts a component of
$\Supp\fr{F}$. Hence $f_*{\down{F}}$ is non-compact. As above,
$K_{f_*{\down{F}}}+\Diff_{f_*{\down{F}}}$ has a regular
complement. By Lemma~\ref{prodolj} this complement can be extended
on $Z$ near $o$ and by Lemma~\ref{podnyal} there exists a regular
complement of $K_X+{\down{F}}$ near $f^{-1}(o)$.
\par
Finally, assume that $f$ is small. Then there exists a
$(K_X+F)$-flip $X\bir X^+$. On $X^+$ the proper transform
$K_{X^+}+{\down{F}}^+$ of $K_X+{\down{F}}$ is antiample or
numerically trivial over $Z$. Again by Lemma~\ref{prodolj} there
is a regular complement on $X^+$ and taking pull-back it on $X$ we
get the desired complement.
\end{proof}
Now we consider the case when $K_X+F$ is not plt.

\begin{lemma}[{\cite{Pr1}, cf.~\cite[9.1]{Sh},
\cite[17.4]{Ut}},~\cite{Sh2}]\label{plt} Let $X$ be a normal
$\QQ$-factorial variety of dimension $\le 3$ and let $F$ be a
boundary on $X$ such that $K_X+F$ is lc, but not plt. Assume that
$F\ne 0$ and $X$ has at worse klt singularities. Then there exists
a plt blow-up $g\colon Y\to X$ such that
\begin{enumerate}
\item
$Y$ is $\QQ$-factorial, $\rho (Y/X)=1$ and the exceptional locus
of $g$ is an irreducible divisor $S$;
\item
$K_Y+S+B=g^*(K_X+F)$ is lc;
\item
$K_Y+S+(1-\ep)B$ is plt and antiample over $X$ for any $\ep>0$.
\end{enumerate}
\end{lemma}


\begin{lemma}
\label{exist-A}
Notation as in \ref{notations}. Let $F$ be such as in
\ref{notation1}. Assume that $K_X+F$ is not plt. Then there exists
a plt blow-up $g\colon (Y,S)\to X$ and a complement $K_Y+S+A$
which is plt.
\end{lemma}
Note that complements $K_Y+S+A$ and $K_X+g_*A$ are exceptional in
the sense of \ref{def-sing}.

\begin{proof}
Apply Lemma~\ref{plt} to $(X,F)$. We get a plt blow-up $g\colon
Y\to X$ with log canonical $g^*(K_X+F)=K_Y+S+B$. Since
$\rho(Y/Z)=2$, the cone $\ov{\NE}(Y/Z)$ is generated by two
extremal rays. Denote them by $R$ and $Q$. One of them, say $R$,
determines the contraction $g\colon Y\to X$ and therefore it is
trivial with respect to $K_Y+S+B$, so $B$ must be positive on $R$.
Since $K_Y+S+B=g^*(K_X+F)$, where $K_X+F$ antiample over $Z$,
$K_Y+S+B$ cannot be nef. Therefore $K_Y+S+B$ is negative on $Q$.
Take small $\beta>0$ and denote $B^\beta:=(1-\beta)B$.
\par
Then $K_Y+S+B^\beta$ is plt and antiample over $Z$. For
sufficiently big and divisible $n\in\NN$ the divisor
$-n(K_Y+S+B^\beta)$ is very ample over $Z$. Take a general member
$L\in |-n(K_Y+S+B^\beta)|$ and denote $A:=B^\beta+(1/n)L$. Then
$K_Y+S+A$ is an $n$-complement of $K_Y+S$. Moreover, $K_Y+S+A$ is
plt, see \cite[4.7]{Ko}.
\end{proof}

\begin{remark}\label{equivalence}
Since $K_Y+S+B^\beta$ is plt and antiample over $Z$, $R^1(g\circ
f)_*\OOO_Y=0$. This gives us that the numerical equivalence of
(Cartier) divisors on $Y$ coincides with $\QQ$-linear one.
\end{remark}

\subsection{}\label{construction-Y}
Now let $f\colon X\to Z$ be a contraction as in \ref{notations},
let $g\colon (Y,S)\to X$ be a plt blow-up from Lemma~\ref{exist-A}
and let $K_Y+S+A$ be a complement on $Y$ which is plt.

Note that $-(K_Y+S)$ is not necessarily nef over $Z$ (so we cannot
apply Lemma~\ref{prodolj} directly). To improve the situation we
consider the following construction.

\begin{proposition}
\label{construction-Y1}
Notation as in~\ref{construction-Y}. One of the following holds.
\begin{enumerate}
\item[\textrm{(A)}]
After a sequence of $(K_Y+S+(1+\ep)A)$-flips $Y\bir\ov{Y}$ over
$Z\ni o$, we get a log variety $(\ov{Y},\ov{S}+(1+\ep)\ov{A})$ and
the diagram
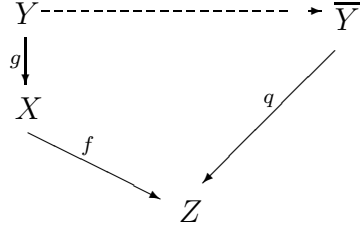
\begin{center}
\begin{picture}(140,100)
\put(10,90){$Y$} \multiput(20,95)(5,0){19}{\line(1,0){3}}
\put(121,95){\vector(1,0){5}} \put(131,88){$\ov{Y}$}
\put(14,85){\vector(0,-1){18}} \put(131,80){\vector(-1,-1){50}}
\put(104,60){$\scriptstyle{q}$} \put(8,75){$\scriptstyle{g}$}
\put(10,54){$X$} \put(15,49){\vector(2,-1){50}}
\put(35,42){$\scriptstyle{f}$} \put(72,15){$Z$}
\end{picture}
\end{center}
such that $K_{\ov{Y}}+\ov{S}$ is plt and $-(K_{\ov{Y}}+\ov{S})$ is
$q$-nef and $q$-big.
\item[\textrm{(B)}]
After a sequence of $(K_Y+S+(1+\ep)A)$-flips $Y\bir\widehat{Y}$
over $Z\ni o$, we get a log variety
$(\widehat{Y},\widehat{S}+(1+\ep)\widehat{A})$ and the diagram
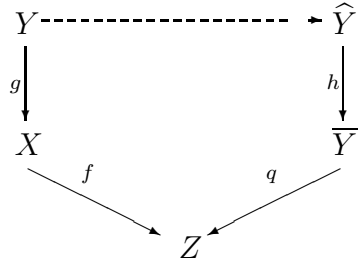
\begin{center}
\begin{picture}(140,100)
\put(10,90){$Y$} \multiput(20,95)(5,0){19}{\line(1,0){3}}
\put(121,95){\vector(1,0){5}} \put(130,90){$\widehat{Y}$}
\put(14,85){\vector(0,-1){28}} \put(8,69){$\scriptstyle{g}$}
\put(134,85){\vector(0,-1){28}} \put(128,69){$\scriptstyle{h}$}
\put(10,44){$X$} \put(130,44){$\ov{Y}$}
\put(15,39){\vector(2,-1){50}} \put(35,35){$\scriptstyle{f}$}
\put(105,35){$\scriptstyle{q}$} \put(133,39){\vector(-2,-1){50}}
\put(72,5){$Z$}
\end{picture}
\end{center}
where $h\colon\widehat{Y}\to\ov{X}$ is a divisorial contraction
which is positive with respect to
$K_{\widehat{Y}}+\widehat{S}\equiv-\widehat{A}$ and negative with
respect to $K_{\widehat{Y}}+\widehat{S}+(1+\ep)\widehat{A}$. In
this case $\ov{S}:=h(\widehat{S})$ is a surface,
$K_{\ov{Y}}+\ov{S}$ is plt and $-(K_{\ov{Y}}+\ov{S})$ is ample
over $Z$.
\end{enumerate}
\end{proposition}

\begin{proof}
For sufficiently small $\ep>0$ the divisor $K_Y+S+(1+\ep)A$ is
also plt (see~\cite[1.3.4]{Sh}, \cite[2.17.4]{Ut}). Apply
$(K_Y+S+(1+\ep)A)$-MMP to $Y$. After a number of flips we get one
of the following:
\begin{itemize}
\item
A log pair $(\ov{Y}, \ov{S}+(1+\ep)\ov{A})$, where
$K_{\ov{Y}}+\ov{S}+(1+\ep)\ov{A})$ is nef over $Z$. Since
$K_{\ov{Y}}+\ov{S}+\ov{A}\equiv 0$, so are both $\ov{A}$ and
$-(K_{\ov{Y}}+\ov{S})$. This is case \textrm{(A)} of
\ref{construction-Y1}.
\item
A log pair $(\widehat{Y}, \widehat{S}+(1+\ep)\widehat{A})$ with a
non-flipping extremal contraction $h\colon\widehat{Y}\to\ov{Y}$
which is negative with respect to
$K_{\widehat{Y}}+\widehat{S}+(1+\ep)\widehat{A}$. This contraction
must be negative with respect to $\widehat{A}$, because
$K_{\widehat{Y}}+\widehat{S}+\widehat{A}\equiv 0$ over $Z$.
Therefore it is divisorial and contracts a component of
$\Supp(\widehat{A})$. In particular, $\ov{S}:=h(\widehat{S})$ is a
divisor. Put $\ov{A}:=\widehat{A}$. Then $-(K_{\ov{Y}}+\ov{S})\qq
\ov{A})$ is big over $Z$. By Kodaira's Lemma (see e.~g.
\cite[0-3-4]{KMM}),
\begin{equation*}
-(K_{\ov{Y}}+\ov{S})\qq (q\text{-ample})+(\text{effective}).
\end{equation*}
If $q$ is not birational, then $-(K_{\ov{Y}}+\ov{S})$ is $q$-ample
because $\rho(\ov{Y}/Z)=1$. If $q$ is birational, then it
contracts $\ov{S}$. This gives us that $f\circ g$ contracts two
divisors. Hence $f$ is a divisorial contraction and $Z\ni o$ is
klt and $\QQ$-factorial \cite[5-1-6]{KMM}. Then
$-(K_{\ov{Y}}+\ov{S})$ must be $q$-ample. We get case \textrm{(B)}
of \ref{construction-Y1}.
\end{itemize}
\end{proof}

\begin{remark}\label{equivalence1}
By \ref{equivalence} the numerical equivalence on $Y$ coincides
with $\QQ$-linear one. All flips and the divisorial contraction
$\widehat{Y}\bir\ov{Y}$ preserve this property, so the same holds
on $\ov{Y}$.
\end{remark}

\begin{proposition}
\label{g(S)-curve}
In notation~\ref{notations} assume that $f(g(S))$ is a curve. Then
$K_X$ has a regular non-exceptional complement.
\end{proposition}
\begin{proof}
Consider case \textrm{(A)} of Proposition~\ref{construction-Y1}.
It is easy to see that $q(\ov{S})=f(g(S))$ is also a curve. Denote
this curve by $C$. We have a contraction $q\colon\ov{S}\to C\ni
o$. It is clear that $-(K_{\ov{S}}+\Diff_{\ov{S}})=
-(K_{\ov{Y}}+\ov{S})|_{\ov{S}}$ is nef over $C$.
By~\cite[3.1]{Sh1} there exists a regular complement of
$K_{\ov{S}}+\Diff_{\ov{S}}$ near $q^{-1}(o)$. By
Lemma~\ref{prodolj} this complement can be extended to a regular
complement $K_{\ov{Y}}+\ov{S}+\ov{D}$ on $\ov{Y}$. The sequence of
maps
\begin{equation*}
\chi\colon\ov{Y}\bir\cdots\bir Y
\end{equation*}
is a sequence of flops with respect to $K_{\ov{Y}}+\ov{S}+\ov{D}$.
They preserve the lc property of $K_{\ov{Y}}+\ov{S}+\ov{D}$
(see~\cite[2.28]{Ut}). This gives us a regular complement
$K_Y+S+D=\chi_*(K_{\ov{Y}}+\ov{S}+\ov{D})$. Finally, $K_X+g_*D$ is
a regular complement on $X$ and $a(S,g_*D)=-1$.
\par
Case \textrm{(B)} of Proposition~\ref{construction-Y1} can be
treated by the similar way: we can construct a complement on
$\ov{S}$, extend it on $\ov{Y}$, pull-back it on $\widehat{Y}$ by
Lemma~\ref{podnyal} and take the proper transform on $X$.
\end{proof}

Thus if $f(g(S))$ is a curve, then we have case (i) of Theorem~A.
From now on we assume that $f(g(S))$ is a point, so both $S$ and
$\ov{S}$ are compact surfaces.

\begin{lemma}
\label{nef-big}
Notation as in~\ref{construction-Y}.  Assume that $f(g(S))$ is a
point. Then $-(K_{\ov{S}}+\Diff_{\ov{S}})$ is nef and big on
$\ov{S}$ (by definition, $(\ov{S},\Diff_{\ov{S}})$ is a weak log
del Pezzo).
\end{lemma}

\begin{proof}
Since $K_{\ov{S}}+\Diff_{\ov{S}}= (K_{\ov{Y}}+\ov{S})|_{\ov{S}}$,
the divisor $-(K_{\ov{S}}+\Diff_{\ov{S}})$ is nef. If
$-(K_{\ov{Y}}+\ov{S})$ is ample over $Z$, then obviously, so is
$-(K_{\ov{S}}+\Diff_{\ov{S}})$. Thus we assume that
$-(K_{\ov{Y}}+\ov{S})$ is not ample over $Z$. This is possible
only in case (A) of \ref{construction-Y1}.
\par
By of Proposition~\ref{construction-Y1} $-(K_{\ov{Y}}+\ov{S})$ is
nef and big. Recall that $\rho(\ov{Y}/Z)=2$, hence $\NE(\ov{Y}/Z)$
is generated by two extremal rays. If $-(K_{\ov{Y}}+\ov{S})$ is
not ample, then it is trivial on some extremal ray $\ov{Q}$ on
$\ov{Y}$. We use the notation of the proof of Lemma~\ref{exist-A}.
Let $\ov{B^\beta}$ be the proper transform of $B^\beta$.
\par
First we consider the case $\ov{B^\beta}\cdot \ov{Q}<0$. We claim
that $\ov{Q}$ is contractible over $Z$. Indeed, $\ov{Q}$ is
negative with respect to $K_{\ov{Y}}+\ov{S}+\ep \ov{B^\beta}$,
where $\ep>0$. Thus $\ov{Q}$ is contractible by Contraction
Theorem~\cite[3-2-1]{KMM}. Let $\ov{g}\colon\ov{Y}\to W$ be the
contraction of $\ov{Q}$ over $Z$. Then
$-(K_{\ov{Y}}+\ov{S})=\ov{g}^*(L)$, where the $\QQ$-divisor $L$ is
ample on $W$, because $\rho (W/Z)=1$. Therefore the linear system
$|-n(K_{\ov{Y}}+\ov{S})|$ is base point free over $Z$ for $n\gg
0$. Further, $\ov{B^\beta}\cdot \ov{Q}<0$, hence the contraction
is birational and contracts a component of $\Supp(\ov{B^\beta})$
or it is small. In both cases the morphism determined by the
linear system $|-n(K_{\ov{S}}+\Diff_{\ov{S}})|$ can contract only
a finite number of curves in $\ov{S}\cap\Supp(\ov{B^\beta})$,
i.~e. $-(K_{\ov{S}}+\Diff_{\ov{S}})$ is big.
\par
Similarly, in the case $\ov{B^\beta}\cdot \ov{Q}>0$ we have that
$\ov{Q}$ is negative with respect to $K_{\ov{Y}}+\ov{S}+\ep
(\ov{A}-\ov{B^\beta})$ (recall that $\ov{A}-\ov{B^\beta}$ is
effective).
\par
Finally, if $\ov{B^\beta}\cdot \ov{Q}=0$, then $K_{\ov{Y}}+\ov{S}$
and $\ov{B^\beta}$ are numerically proportional, because
$\rho(\ov{Y}/Z)=2$. By \ref{equivalence1} $K_{\ov{Y}}+\ov{S}\qq
c\ov{B^\beta}$ for some $c\in\QQ$. Whence $K_Y+S\qq cB^\beta=
c(1-\beta)B$ and $K_Y+S+B\equiv (1+c(1-\beta))B$. Since $-(K_Y+S)$
is $g$-ample and $K_Y+S+B$ is $g$-numerically trivial, we have
$1+c(1-\beta)=0$. So $g^*(K_X+F)=K_Y+S+B\equiv 0$, a
contradiction. This proves the lemma.
\end{proof}

\begin{proposition}
\label{m-n}
Let $f\colon X\to Z$ be as in~\ref{notations}. Then the canonical
divisor $K_X$ is $n$-complementary for some $n<\NNN_2$, where
$\NNN_2$ is the constant from Theorem~\ref{Const}.
\end{proposition}
\begin{proof}
First we consider case \textrm{(A)} of
Proposition~\ref{construction-Y1}. By Shokurov's
theorem~\ref{Const} we have an $n$-complement
$K_{\ov{S}}+\ov{\Theta}$ of $K_{\ov{S}}+\Diff_{\ov{S}}$, where
$n<\NNN_2$. By Lemma~\ref{prodolj} $K_{\ov{S}}+\ov{\Theta}$ can be
extended to an $n$-complement $K_{\ov{Y}}+\ov{S}+\ov{D}$ on
$\ov{Y}$. As in the proof of Proposition~\ref{g(S)-curve} its
proper transform on $X$ is an $n$-complement of $K_X$.
\par
In case \textrm{(B)} of Proposition~\ref{construction-Y1} we have
a contraction $q\colon\ov{X}\to Z$ and a boundary
$\ov{S}:=h(\widehat{S})$ such as in Lemma~\ref{semistable}.
Therefore $K_{\ov{X}}+\ov{S}$ is $n$-complementary for $n<\NNN_2$.
We can pull-back this complement on $\ov{Y}$ by
Lemma~\ref{podnyal}. The rest is similar to case \textrm{(A)}.
\end{proof}

\begin{proposition}
\label{non-exceptional}
Notation as above. Assume that there exists a complement
$K_{\ov{S}}+\ov{\Theta}$ of $K_{\ov{S}}+\Diff_{\ov{S}}$ which is
not klt. Then $K_X$ has a regular non-exceptional complement.
\end{proposition}
\begin{proof}
Assume that $K_{\ov{S}}+\ov{\Theta}$ is not klt. As in the proof
of Proposition~\ref{g(S)-curve} it is sufficient to show that
$K_{\ov{S}}+\Diff_{\ov{S}}$ has a regular non-klt complement.

The following is an easy consequence of~\cite[2.3]{Sh1}
(see~\cite[Corollary~11]{Pr1}).

\begin{lemma}
\label{qua}
Let $S$ be a projective normal surface and let $\Theta\geq\Delta$
be boundaries on $S$ such that
\begin{enumerate}
\item
coefficients of $\Delta$ are standard;
\item
$K_S+\Delta$ is klt;
\item
$-(K_S+\Delta)$ is ample;
\item
$K_S+\Theta$ is lc but not klt;
\item
$-(K_S+\Theta)$ is nef.
\end{enumerate}
Then $K_S+\Delta$ has a regular complement which is not klt.
\end{lemma}

By this lemma $K_{\ov{S}}+\Diff_{\ov{S}}$ has a regular non-klt
complement. By Lemma~\ref{prodolj} it can be extended on $\ov{X}$.
The proper transform of it on $X$ gives us a regular complement
which is non-exceptional.
\end{proof}

We obtain very important
\begin{corollary}
\label{klt1}
Let $f\colon X\to Z$ be as in~\ref{notations}. Let $g\colon
(Y,S)\to X$ be a plt blow-up constructed in \ref{exist-A} and let
$K_Y+S+D$ be some complement. If $K_Y+S+D$ is not plt, then
$K_Y+S$ has a regular non-plt complement.
\end{corollary}
\begin{proof}
Notation as in the proof of Proposition~\ref{construction-Y1}. We
claim that the proper transform $K_{\ov{Y}}+\ov{S}+\ov{D}$ of
$K_Y+S+D$ on $\ov{Y}$ is not plt. Indeed, in case (A) $Y\to\ov{Y}$
is a sequence of flops with respect to $K+S+D$ and the fact
follows by~\cite[2.28]{Ut}. Similarly, in case (B)
$K_{\widehat{Y}}+\widehat{S}+\widehat{D}$ is not plt. The
contraction $h\colon\widehat{Y}\to\ov{X}$ is log crepant with
respect to $K_{\widehat{Y}}+\widehat{S}+\widehat{D}$, so
$K_{\ov{Y}}+\ov{S}+\ov{D}$ is not plt by \cite[3.10]{Ko}.
\par
By Lemma~\ref{connect-(2)3} below $K_{\ov{Y}}+\ov{S}+\ov{D}$ is
not plt near $\ov{S}$. Therefore $K_{\ov{S}}+
\Diff_{\ov{S}}(\ov{D})$ is not klt. Now the corollary follows by
Proposition~\ref{non-exceptional}.
\end{proof}

\begin{lemma}\label{connect-(2)3}
Let $f\colon X\to Z$ be a contraction from a threefold with $\dim
Z>0$ and let $D$ be a boundary on $X$. Fix a point $o\in Z$.
Assume that
\begin{enumerate}
\item
$K_X+D$ is lc and not plt near $f^{-1}(o)$;
\item
$K_X+D\equiv 0$ over $Z$;
\item
there is an irreducible component $S_0\subset\down{D}$ such that
$f(S_0)\ne Z$.
\end{enumerate}
Then $K_X+D$ is not plt near $S_0\cap f^{-1}(o)$.
\end{lemma}

In the proof we follow the proof of~\cite[6.9]{Sh}\footnote{The
referee pointed out that Lemma~\ref{connect-(2)3} easily follows
from \cite[Proposition~2.1]{F}}.

\begin{proof}
If $f$ is birational, then the lemma follows by
Theorem~\ref{connect}. We assume that $Z$ is one or
two-dimensional. Consider $Z$ as a sufficiently small neighborhood
of $o$. Replace $(X,D)$ with its $\QQ$-factorial dlt model (see
\cite[9.1]{Sh} or \cite[8.2.2]{Ut}). Put $B:=\fr{D}$ and
$S:=\down{D}$. Since $K_X+S+B$ is not plt, $S$ has at least two
irreducible components. Assume that $K_X+S+B$ is plt near $S_0$.
Then there are no irreducible components $S_i\ne S_0$ of $S$
meeting $S_0$. Hence $S_0$ is a connected component of $S$. Thus
we have
\begin{enumerate}
\item[(i)${}'$]
$K_X+S+B$ is lc over the germ $o\in Z$ and $X$ is $\QQ$-factorial;
\item[(ii)${}'$]
$K_X+S+B\equiv 0$ over $Z$;
\item[(iii)${}'$]
$S_0$ is a connected component of $S$, it is irreducible and
$f(S_0)\ne Z$;
\item[(iv)${}'$]
$LCS(X,S+B)=S$ and $S\ne S_0$.
\end{enumerate}
If $K_X+B\equiv -S$ is nef, then  $S\supset f^{-1}(o)$. Indeed,
otherwise there is a curve $C\subset f^{-1}(o)$ such that $C\cdot
S>0$. Obviously, in this case $S$ is connected near $f^{-1}(o)$, a
contradiction. Now we assume that $K_X+B$ is not nef. Run
$(K_X+B)$-MMP over $Z$. On each step $K_X+B$ is klt. Hence
$K_X+S+B$ is klt outside of $S$ and $LCS(X,S+B)=S$. Let $R$ be a
$(K_X+B)$-negative extremal ray and let $\psi\colon X\to X'$ be
its contraction. Since $R\cdot S>0$, $\psi$ cannot contract a
connected component of $S$. In particular, $\psi$ cannot contract
$S_0$. Further, by Theorem~\ref{connect} $S$ is connected near
each fiber of $\psi$. Therefore the number of connected components
of $S$ remains the same. We have shown that properties
(i)${}'$--(iv)${}'$ hold on every step of MMP.

\par
At the end we get a non-birational $(K_X+B)$-negative extremal
contraction $f'\colon X\to Z'$ over $Z$. Then $S\equiv -(K_X+B)$
is $f'$-ample. On the other hand, $f'(S_0)\ne Z'$, so $S_0\equiv
0$ over $Z'$. Write $S=\sum S_i+\sum S_j$, where $S_i$'s are
$f'$-ample and $S_j$'s are $f'$-numerically trivial components. By
the above remarks, both terms are non-empty. Any component of type
$S_i$ intersects all curves in fibers of $f'$. In particular, it
intersects any component of type $S_j$. Therefore $S$ is
connected. This contradicts to our assumptions
(iii)${}'$--(iv)${}'$.
\end{proof}

\begin{definition}
Let $f\colon X\to Z\ni o$ be a contraction such as in
\ref{notations}. Then it is said to be \textit{exceptional} if any
non-klt complement $K_X+D$ near $f^{-1}(o)$ is exceptional.
\end{definition}

\begin{corollary}
\label{non-klt-exceptional}
Let $f\colon X\to Z$ be as in~\ref{notations} and assume that it
is not exceptional. Then $K_X$ has a regular complement which is
also non-exceptional.
\end{corollary}

\section{Proof of Theorem~A}\label{proof}
In this section we prove Theorem~A.

\begin{proposition}
\label{deg<}
Notation as in~\ref{notations} and \ref{notation1}. Assume that
$K_X+F$ is plt. Let $\psi\colon X'\to X$ be a finite \'etale in
codimension one cover. Assume that all complements of
$K_X+\down{F}$ are exceptional (in this case it means that they
are plt). Then $\deg\psi$ is bounded by an absolute constant.
\end{proposition}
In case when $K_X+F$ is plt this assertion gives us a stronger
result than (ii) of Theorem~A because here the absolute constant
does not depend on $\ep$.

\begin{proof}
By Lemma~\ref{semistable} we may assume that $\down{F}$ is compact
and $-(K_X+{\down{F}})$ is $f$-ample. Thus ${\down{F}}=f^{-1}(o)$.
Denote $\Delta :=\Diff_{\down{F}}$, $m:=\deg\psi$. If
$K_{\down{F}}+\Delta$ has a regular complement, then by
Lemma~\ref{prodolj} it can be extended on $X$, so we assume the
opposite. By Theorem~\ref{Const} $K_{\down{F}}+\Delta$ has an
$n$-complement $K_{\down{F}}+\Theta$ for $n<\NNN_2$. Moreover, we
may assume that $K_{\down{F}}+\Theta$ is klt (see
Lemma~\ref{qua}). Therefore both $K_{\down{F}}+\Theta$ and
$K_{\down{F}}+\Delta$ are $(1/\NNN_2)$-lt, because
$n(K_{\down{F}}+\Theta)\sim 0$. Now put $F':=\psi^{-1}F$. It is
easy to see $\down{F'}=\psi^{-1}(\down{F})$. Since $\psi\colon
X'\to X$ is \'etale in codimension one,
\begin{equation*} K_{X'}+F'=\psi^*(K_X+F), \qquad
K_{X'}+\down{F'}=\psi^*(K_X+{\down{F}}). \end{equation*} By
\cite[\S 2]{Sh} (see also \cite[20.3]{Ut}, \cite[3.16]{Ko}) we
have that $K_{X'}+F'$ is plt. So, by Inversion of Adjunction and
by Connectedness Lemma $\down{F'}$ is irreducible and normal. It
is clear that the restriction $\var=\psi|_{\down{F'}}\colon
\down{F'}\to {\down{F}}$ is a finite morphism of the same degree
as $\psi$. Denote $\Delta':=\Diff_{\down{F'}}$. Thus by
Adjunction~\ref{Inv-Adj} we have
$K_{\down{F'}}+\Delta'=\var^*(K_{\down{F}}+\Delta)$. Therefore
$(\down{F'}, \Delta')$ is also $(1/\NNN_2)$-lt (see~\cite[Sect.
2]{Sh} and also~\cite[20.3]{Ut}). Similarly,
$K_{\down{F'}}+\Theta'=\var^*(K_{\down{F}}+\Theta)$, where
$K_{\down{F'}}+\Theta'$ is an $n$-complement of
$K_{\down{F'}}+\Delta'$. By~\cite{A2} the surface $\down{F'}$
belongs to a finite number of families. Since $\Theta'\in
|-nK_{\down{F'}}|$, the set $\{(\down{F'}, \Theta')\}$ is also
bounded. Taking into account $\Delta'\leq \Theta'$, we have that
so is $\{(\down{F'}, \Delta')\}$. Therefore
$0<(K_{\down{F'}}+\Delta')^2<\Const$. Finally, the equality
$(K_{\down{F'}}+\Delta')^2=(\deg\var)(K_{\down{F}}+\Delta)^2$
gives us that $\deg\var=\deg\psi$ is bounded and proves the
proposition.
\end{proof}

The following is the main step in the proof of Theorem~A.

\begin{proposition}
\label{claim}
Let $f\colon X\to Z\ni o$ be an contraction as in~\ref{notations}.
Assume that $X$ has only $\ep$-lt singularities, where $\ep>0$.
Assume also that $f\colon X\to Z\ni o$ is exceptional. Let
$\psi\colon X'\to X$ be a finite Galois \'etale in codimension one
cover with connected $X'$. Then the degree of $\psi\colon X'\to X$
is bounded by a constant depending only on $\ep$.
\end{proposition}

\begin{proof}
First note that since $X$ is a germ along $f^{-1}(o)$,
$(f\circ\psi)^{-1}(o)$ is connected.
\par
Below we will use notation of \ref{construction-Y}. Taking into
account Proposition~\ref{deg<} we may assume that $K_X+F$ is not
plt. Let $g\colon (Y,S)\to X$ be a plt blow-up from
Lemma~\ref{exist-A}. By Proposition~\ref{g(S)-curve} we may assume
that $f\circ g(S)$ is a point and By Proposition~\ref{m-n} there
exists an $n$-complement $K_Y+S+D$, where $n<\NNN_2$. Denote
$\Delta :=\Diff_S$ and $\Theta :=\Diff_S(D)$. By our assumptions
and by \ref{klt1} $K_Y+S+D$ is plt hence $K_S+\Theta$ is klt.
Since $n(K_S+\Theta)\sim 0$, $K_S+\Theta$ is also $(1/\NNN_2)$-lt.
As in the proof of Proposition~\ref{deg<} by~\cite{A2} $S$ lies in
a finite number of families and because $n\Theta\in |-nK_S|$, the
pair $(S,\Theta)$ also lies in a finite number of families. Since
$\Delta\le\Theta$ and $\Delta =\sum (1-1/m_i)\Delta_i$ has only
standard coefficients, so is $(S,\Delta)$. Therefore
\subsection{We may assume that $(S,\Delta)$ is fixed.}
\label{fixed1}
Let $Y'$ be the normalization of a dominant component of
$Y\times_XX'$. Consider the commutative diagram
\begin{equation}
\label{CD1}
\begin{CD}
Y'@>\phi>>Y\\ @Vg'VV @VgVV\\ X'@>\psi>>X\\ @Vf'VV @VfVV\\
Z'@>\pi>>Z\\
\end{CD}
\end{equation}
where $X'\to Z'\to Z$ is the Stein factorization. It is clear that
$\pi^{-1}(o)$ is a point (because $(f\circ\psi)^{-1}(o)$ is
connected). We claim that $g'\colon Y'\to X'$ is a plt blow-up. It
is clear that $\phi\colon Y'\to Y$ is a finite Galois morphism and
its ramification divisor can be supported only in $S$. Denote
$S':=\phi^{-1}(S)$. Then $S'$ is the exceptional divisor of the
blow-up $g'\colon Y'\to X'$. We have that
\begin{equation}
\label{for-2}
K_{Y'}+S'=\phi^*(K_Y+S)
\end{equation}
and this divisor is plt \cite[2.2]{Sh}, \cite[20.3]{Ut}. This
gives us that $S'$ is normal (see \ref{Inv-Adj}). On the other
hand, $-(K_{Y'}+S'+B')=-\phi^*(K_Y+S+B)$ if nef and big over $Z'$.
By Connectedness Lemma~\ref{connect} $S'$ is connected near each
fiber of $f'\circ g'$ (and hence it is irreducible). Thus we have
that $g'\colon Y'\to X'$ is a plt blow-up. This shows also that
$K_{X'}+f'_*D'$ is exceptional. Similar to (\ref{for-2}) we have
\begin{equation}
\label{for-3}
K_{Y'}+S'+D'=\phi^*(K_Y+S+D)
\end{equation}
and again this divisor is plt near $S'$ \cite[2.2]{Sh},
\cite[20.3]{Ut}. Put $\Theta':=\Diff_{S'}(\Theta)$ and
$\Delta':=\Diff_{S'}$. Thus we have a finite Galois morphism
$\var=\phi |_{S'}\colon S'\to S$ such that
\begin{equation}
\label{two-eq}
K_{S'}+\Delta'=\var^*(K_S+\Delta)\quad\text{ and}\quad
K_{S'}+\Theta'=\var^*(K_S+\Theta).
\end{equation}

\subsection{We may assume that $(S',\Delta')$ is fixed.}
\label{fixed2}
Indeed, by~\cite[Sect. 2]{Sh} (see also~\cite[20.3]{Ut}) taking
into account (\ref{two-eq}) we have $discr(S', \Theta')\geq
discr(S, \Theta)$. As in \ref{fixed1} by~\cite{A2} the set
$\{(S',\Delta')\}$ is bounded.

\begin{lemma}
\label{deg_var_is_bounded}
$\deg\var$ is bounded.
\end{lemma}

\begin{proof}
First we assume that $g(S)$ is a point. In this case both
$(S,\Delta)$ and $(S',\Delta')$ are log del Pezzo surfaces. We
have $0<(K_{S'}+\Delta')^2=(\deg\var)(K_S+\Delta)^2$. Since for
$(K_{S'}+\Delta')^2$ and $(K_S+\Delta)^2$ there are only a finite
number of possibilities, $\deg\var $ is bounded. In the case when
$g(S)$ is a curve $(S, \Xi)$ is a log del Pezzo for $\Xi
=\Diff_S(B^\beta)\geq\Delta$ (see the proof of
Lemma~\ref{exist-A}).

\begin{lemma}
\label{max-DelPezzo}
Let $(S,\Delta)$ be a projective log surface. Consider the set of
boundaries on $S$

\begin{equation*}
\MMM(S,\Delta):=\{\Xi : \Xi\ge\Delta,\ K_S+\Xi\ \text{is klt\
and}\ -(K_S+\Xi) \text{is nef}\}.
\end{equation*}
Assume that $\MMM(S,\Delta)\ne\emptyset$. Then there exists a
boundary $\Xi_{\max}$ such that $\Xi_{\max}\ge\Delta$,
$-(K_S+\Xi_{\max})$ is nef, $K_S+\Xi_{\max}$ is lc, and
$(K_S+\Xi_{\max})^2\ge(K_S+\Xi)^2$ for any $\Xi\in\MMM(S,\Delta)$.
\end{lemma}

Note that $\Xi_{\max}$ is contained in the closure
$\ov{\MMM(S,\Delta)}$.

\begin{proof}
We may assume that $-(K_S+\Xi)$ is big, otherwise $(K_S+\Xi)^2=0$
for all $\Xi\in\MMM(S,\Delta)$. Let
$\Xi=\sum\xi_j\Xi_j\in\MMM(S,\Delta)$ and let $\Xi_i$ be a
component of $\Xi$ such that $\Xi_i^2\ge 0$. Then we have
$\Xi-\alpha\Xi_i\in\MMM(S,\Delta)$ if $\Xi-\alpha\Xi_i\ge\Delta$.
Indeed, if $(K_S+\Xi-\alpha\Xi_i)\cdot L>0$ for some irreducible
curve $L$, then $\Xi_i\cdot L<0$, a contradiction. Further,
\begin{multline*}
(K_S+\Xi-\alpha\Xi_i)^2=\\ (K_S+\Xi)^2-\alpha (K_S+\Xi)\cdot\Xi_i-
\alpha (K_S+\Xi-\alpha\Xi_i)\cdot\Xi_i\ge (K_S+\Xi)^2.
\end{multline*}
So we may assume that for any component $\Xi_i$ we have $\xi_i=0$
whenever $\Xi_i^2\ge 0$ and $\Xi_i\not\subset\Supp(\Delta)$. In
particular, all components of $\Xi$ are contained in
$\Supp(\Delta)$ or have negative self-intersection number. It is
easy to show that $\NE(S)$ is polyhedral (see
e.~g.~\cite[2.5]{Sh1}, \cite[4.12]{Pr1}). Hence there is only a
finite number of curves with negative self-intersections. Thus we
may assume that components of $\Xi$ belong to a finite set of
curves. Finally, we find the maximum of the quadratic form
$(K_S+\Xi)^2$ on the compact set
$\ov{\MMM(S,\Delta)}\cap(\oplus\RR{\cdot}[\Xi_i])$. Obviously,
this maximal $\Xi_{max}$ is rational.
\end{proof}
To finish the proof of Lemma~\ref{deg_var_is_bounded} we take a
boundary $\Xi_{\max}$ on $S$ such as in \ref{max-DelPezzo} and
define a boundary $\Xi'_m$ on $S'$ by the formula
$K_{S'}+\Xi'_m=\var^*(K_S+\Xi_{\max})$ (see~\cite[2.1]{Sh},
\cite[20.2]{Ut}). It is clear that $\Xi'_m\ge\Delta'$ (because
$\Xi_{\max}\ge\Delta$). Again we have
\begin{equation*}
\Const_2\ge(K_{S'}+\Xi'_m)^2=
(\deg\var)(K_S+\Xi_{\max})^2\ge(\deg\var) \Const_1>0
\end{equation*}
This gives us $\deg\var\le \Const_2/\Const_1$ and proves
Lemma~\ref{deg_var_is_bounded}.
\end{proof}

Now we finish the proof of Proposition~\ref{claim}. Let $G$ be the
Galois group of $\psi\colon X'\to X$. Then $G$ acts on $Y'$ and
$Y=Y'/G$. Moreover, $S'$ is $G$-invariant, so we have the action
of $G$ on $S'$. Let $G_0\subset G$ be a subgroup which acts on
$S'$ trivially. It is clear that $G_0\subset G$ is a normal
subgroup of index $[G:G_0]=\deg\var$. Let $G_1:=G/G_0$. Then
$S=S'/G_1$. We have proved that the order of $G_1$ is bounded. It
remains to show that the order of $G_0$ is also bounded. Denote it
by $r$. Let us write
\begin{equation}
\label{two-eq2}
K_Y=g^*K_X+aS, \qquad K_{Y'}=g^{\prime *}K_{X'}+a'S'.
\end{equation}
By assumption $a>-1+\ep$. On the other hand, $\phi^*S=rS'$ and by
the Hurwitz formula
\begin{equation}
\label{two-eq3}
K_{Y'}=\phi^*K_Y+(r-1)S',\qquad K_{X'}=\psi^*K_X.
\end{equation}
Combining (\ref{two-eq2}) and (\ref{two-eq3}) it is easy to obtain
(cf.~\cite[\S 2]{Sh}, \cite[3.16]{Ko})
\begin{equation}
\label{one-eq4}
a'+1=r(a+1)>r\ep
\end{equation}
Let us consider a sufficiently general curve $H$ on $S'$. If
$g'(S')$ is a point, then we can take $H\in
|-n_0(K_{S'}+\Delta')|$, where $n_0$ depends only on
$(S',\Delta')$. If $g'(S')$ is a curve, then we can take $H$ as
the general fiber of $S'\to g'(S')$. From (\ref{two-eq2}) we have
\begin{equation*}
\label{(osn-discr0)}
K_{S'}+\Delta'\equiv (K_{Y'}+S')|_{S'}\equiv (1+a')S'|_{S'}.
\end{equation*}
This gives us
\begin{equation}\label{(osn-discr)}
H\cdot (K_{S'}+\Delta')=H\cdot (K_{Y'}+S')=(1+a')H\cdot S',
\end{equation}
where the left side is equal to $-n_0(K_{S'}+\Delta')^2$ if
$g'(S')$ is a point and $-2+H\cdot \Delta'\ge -2$ if $g'(S')$ is a
curve. In particular, $H\cdot (K_{S'}+\Delta')$ depends only on
$(S',\Delta')$, but not on $Y'$. Recall that the coefficients of
$\Delta'=\Diff_{S'}$ are standard (see~\cite[3.9]{Sh},
\cite[16.6]{Ut}), so we can write $\Delta'=\sum_{i=1}^r
(1-1/m_i)\Delta_i'$ where $m_i\in\NN$, $r\ge 0$. Put
$m':=\mt{l.c.m.}(m_1,\dots , m_r)$. Again by~\cite[3.9]{Sh} both
$m'S'$ and $m'(K_{S'}+\Delta')$ are Cartier along $H$. So we can
rewrite (\ref{(osn-discr)}) as $N=(1+a')k$, where $N=-m'H\cdot
(K_{S'}+\Delta')$ is a fixed natural number and $k=-m'(H\cdot S)$
is also natural. Thus by (\ref{one-eq4}) $N=(1+a')k>kr\ep\ge
r\ep$. This gives us that $r<N/\ep$ is bounded and proves the
proposition.
\end{proof}

\begin{remark*}[Suggested by the referee]
On can prove Proposition~\ref{claim} by working on the
normalization of $\ov{Y}$ in the function field of $Y'$. In this
way Lemma~\ref{deg_var_is_bounded} becomes immediate from
\cite{A2} and the case division by $\dim g(S)$ is not needed.
\end{remark*}

\begin{proof}[{Proof of Theorem~A}]
Let $f\colon X\to Z$ be a contraction from Theorem~A. So, $f$ is a
log conic bundle or a log del Pezzo fibration. Assume that $K_X$
has no regular non-klt complements. Then $f$ is exceptional by
\ref{non-klt-exceptional} and it is sufficient to prove that
singularities $(Z\ni o)$ (resp. multiplicities of $f^{-1}(o)$) are
bounded in terms of $\ep$.

First, consider the case $\dim Z=2$. It is known that $(Z\ni o)$
is a quotient singularity, i.~e. $(Z,o)=(Z',o')/G$, where $ (Z'\ni
o')$ is a non-singular germ and $G$ is a finite group acting on
$Z'$ free outside $o$. We may assume that $G\subset GL_2(\CC)$.
Two conjugate subgroups give us analytic isomorphic singularities
$(Z\ni o)$. So it is sufficient to show boundedness of $|G|$. Let
us construct the base change
\begin{equation}
\label{CD-main}
\begin{CD}
X'@>\psi>>X\\ @Vf'VV @VfVV\\ Z'@>\pi>>Z\\
\end{CD}
\end{equation}
where $X'$ is the normalization of the dominant component of
$X\times_ZZ'$. It is clear that $\psi\colon X'\to X$ is a finite
Galois morphism which is \'etale in codimension $1$. Take a
boundary $F$ as in \ref{notation1}. If $K_X+F$ is plt, then the
assertion follows by \ref{deg<}. Moreover, in our case there
exists a regular complement of $K_X$ (because $\down{F}$ is
non-compact). If $K_X+F$ is not plt, then the degree of $\psi$ is
bounded by Proposition~\ref{claim}.
\par
In the case $\dim Z=1$ denote $L:=f^{-1}(o)_{\mt{red}}$. Then
$mL\sim 0$, where $m$ is the multiplicity of $f^{-1}(o)$. This
gives us the $m$-cyclic cover $\psi\colon X'\to X$ which is
\'etale in codimension one. The rest is similar to the case above.
\end{proof}

\section{Examples and concluding remarks}
\label{last}
\begin{proposition}[cf. {\cite[2.7]{MP}}]
\label{unique}
Let $f\colon X\to Z$ be a contraction as in Theorem~A. Assume $f$
is exceptional. Then there exists a divisor $S$ of $K(X)$ such
that $a(S,F)=-1$ for every non-klt complement $K_X+F$.
\end{proposition}

\begin{proof}
Let us take any complement $K_X+F$ which is not klt. Then there is
a (not necessarily exceptional) divisor $S$ of the field $K(X)$
with discrepancy $a(S,F)=-1$. Let us take another non-klt
complement $K_X+F^1$. We claim that $a(S,F^1)=-1$. Indeed, we can
construct a continuous family of boundaries $F(t)$, $t\in [0,1]$
such that $F(0)=F$, $F(1)=F^1$, $K_X+F(t)$ is lc, non-klt and
numerically trivial for all $t\in [0,1]$. To see this take an
effective divisor $L$ on $Z$ containing $o$ and put
$F(t)=(1-t)F+tF^1+c(t)f^*L$, where $c(t)$ is the log canonical
threshold of $(X,(1-t)F+tF^1)$ with respect to $f^*L$
(see~\cite{Sh}, \cite{Ko}). Since $c(t)$ can be computed from a
finite number of linear inequalities (see~\cite[8.5]{Ko}), $c(t)$
is a picewise linear function in $t$. In particular, it is
continuous. Obviously, $K_X+F(t)\equiv 0$. Further, for any
divisor $E$ of $K(X)$ the discrepancy $a(E,F(t))$ is also
continuous in $t$. Put
\begin{equation*}
t_0:=\sup \{ t\in [0,1] : a(S,F(t))=-1\}.
\end{equation*}
By continuity, there is another divisor $E$ such that
$a(E,F(t_0))=-1$. This contradicts to
Corollary~\ref{non-klt-exceptional}.
\end{proof}

Below we consider explicit examples of contractions such as in
Theorem~A and we shall show that both cases (i) and (ii) of
Theorem~A can really occur.

\begin{example}
\label{example1}
Let $(Z\ni o)$ be a two-dimensional klt singularity, let
$X:=\PP^1\times Z$, and let $f\colon X\to Z$ be the projection.
Then $f$ is a generically conic bundle as in Theorem~A. It is
known that $K_Z$ is $1$, $2$, $3$, $4$ or $6$-complementary in
cases $\AA_n$, $\DD_n$, $\EE_6$, $\EE_7$, $\EE_8$,
respectively~\cite[5.2.3]{Sh}. Let $K_Z+\Delta$ be this
$n$-complement with minimal $n$ and let $S$, $S_1$ be two distinct
sections of $f$. Then $K_X+f^*\Delta+S+S_1$ is an $n$-complement
on $X$. This $n$ is also minimal. Indeed, let $K_X+B$ be any
$m$-complement on $X$. Take a general point $P$ on the central
fiber of $f$ and consider the general section $H\ni P$ of $f$.
Then $K_H+\Diff_H(f^*\Delta)$ is an $m$-complement on $H$ near
$P$. Since $(H\ni P)\simeq (Z\ni o)$, we have that $K_Z$ is also
$m$-complementary.
\end{example}

\begin{example}
Let $G\subset GL_2(\CC)$ be a finite subgroup without
quasireflections. Consider the natural projection $f\colon
(\CC^2\times\PP^1)/G\to \CC^2/G$, where the action of $G$ on
$\PP^1$ is induced from $\CC^2$. Then $f$ satisfies conditions of
Theorem~A. As in the example above one can show that there exists
a regular complement. Similarly, the base surface here can have
any lt singularity.
\end{example}

The examples above shows that the condition of $\ep$-lt cannot be
removed from Theorem~A.

In the case $\dim Z=1$ we can construct a huge number of examples
as products $F\times\CC^1\to\CC^1$, where $F$ is a del Pezzo with
klt singularities. If we take $F$ so that $K_F$ has no regular
complements, we get an exceptional contraction.

The case of birational contractions such as in~\ref{notations} was
studied by Shokurov~\cite[Sect. 7]{Sh}. Here we note only that the
same arguments as in Sect.~\ref{proof} allow us to show that for
divisorial contractions $f\colon X\to Z\ni o$, where $X$ has only
$\ep$-lt singularities one has either $K_X$ has a regular
complement or the index of $(Z\ni o)$ is bounded by $\CCC(\ep)$.

Almost all arguments can be applied also to the study of log Fano
varieties (i.~e. the case when $Z$ is a point). The problem here
is to construct a boundary $F$ as in~\ref{notations}. Then we have
to find a very singular element in $|-nK_X|$ for some $n\in\NN$
(so called ``tiger''). Such a divisor can be constructed by using
Riemann-Roch if $-K_X^3$ is sufficiently large. In the terminal
case we can use also Riemann-Roch for Weil divisors
(see~\cite{K1}, \cite{A1}).

\subsection*{Acknowledgements}
I am grateful to Professor V.~V.~Shokurov for useful discussions
and criticisms. I thank Professor S.~Mori for valuable remarks.
Finally, I would also like to thank the referee for pointing out
some mistakes and misprints, especially in Lemmas~\ref{qua},
\ref{connect-(2)3} and Proposition~\ref{claim}.

\end{document}